\makeatletter

\def\eqalign#1{\,\vcenter{\openup\jot\m@th
  \ialign{\strut\hfil$\displaystyle{##}$&$\displaystyle{{}##}$\hfil
      \crcr#1\crcr}}\,}
\def\eqalignno#1{\displ@y \tabskip\@centering
  \halign to\displaywidth{\hfil$\displaystyle{##}$\tabskip\z@skip
    &$\displaystyle{{}##}$\hfil\tabskip\@centering
    &\llap{$##$}\tabskip\z@skip\crcr
    #1\crcr}}
\def\leqalignno#1{\displ@y \tabskip\@centering
  \halign to\displaywidth{\hfil$\displaystyle{##}$\tabskip\z@skip
    &$\displaystyle{{}##}$\hfil\tabskip\@centering
    &\kern-\displaywidth\rlap{$##$}\tabskip\displaywidth\crcr
    #1\crcr}}

\makeatother

\newdimen\pixel \pixel=.00333333 true in

\documentclass[11pt]{article}
\usepackage{epsfig}
\usepackage{fullpage}

\def \bigpar {}
\def \smallpar {}   

\newlength{\saveindent}
\newenvironment{proof}%
      {\bigpar{\bf Proof:}\ 
             \setlength{\saveindent}{\parindent} 
                       \ignorespaces}%
      {\stopproof\ignorespaces\bigbreak \setlength{\parindent}{\saveindent}}

      {\bigpar{\bf Proof:}%
             \setlength{\saveindent}{\parindent} 
                       \,\ignorespaces}%
      {\ignorespaces\bigbreak \setlength{\parindent}{\saveindent}}

      {\bigpar{\bf Proof:}\ %
             \setlength{\saveindent}{\parindent} 
                       \ignorespaces}%
      {\ignorespaces\bigbreak \setlength{\parindent}{\saveindent}}

\newenvironment{proofof}[1]%
      {\bigpar{\bf#1:}\ %
             \setlength{\saveindent}{\parindent} 
                       \ignorespaces}%
      {\stopproof\ignorespaces\bigbreak \setlength{\parindent}{\saveindent}}

\newenvironment{remark}%
      {\smallpar{\bf Remark:}\ 
                       \ignorespaces}%
      {\stopproof\ignorespaces\medbreak \setlength{\parindent}{\saveindent}}

\newenvironment{remark*}%
      {\smallpar{\bf Remark:}\ 
                       \ignorespaces}%
      {\ignorespaces\medbreak \setlength{\parindent}{\saveindent}}

      {\smallpar{\bf Remarks:}\ 
                       \ignorespaces}%
      {\stopproof\ignorespaces\medbreak \setlength{\parindent}{\saveindent}}

\newenvironment{remarks*}%
      {\smallpar{\bf Remarks:}\ 
                       \ignorespaces}%
      {\ignorespaces\medbreak \setlength{\parindent}{\saveindent}}

      {\smallpar{\bf Remark:}\ %
                       \ignorespaces}%
      {\stopproof\ignorespaces\medbreak \setlength{\parindent}{\saveindent}}

      {\smallpar{\bf Remarks:}\ %
                       \ignorespaces}%
      {\stopproof\ignorespaces\medbreak \setlength{\parindent}{\saveindent}}
\makeatother

\newtheorem{theorem}{Theorem}
\newtheorem{lemma}[theorem]{Lemma}

\newtheorem{corollary}[theorem]{Corollary}

\newtheorem{example}{Example}

\def\begex{\begin{example}\parindent=0pt \rm}
\def\endex{\end{example}}
\def\square{\vbox{\hrule height.2pt\hbox{\vrule width.2pt height5pt \kern5pt
                                   \vrule width.2pt} \hrule height.2pt}}
\def\stopproof{\hfill \square \smallskip}

\def \choosek {{ {\textstyle {2^d \choose k}^{-1}}}}
\def \choosekm {{ {\textstyle {2^d \choose k-1}^{-1}}}}
\def \f {{\cal F}}

\def \lam {{\Lambda}}

\def \lone#1{{||#1||_1}}
\def \ltwo#1{{||#1||_2}}
\def \ss#1{{||#1||^2_2}}
\def\sfrac#1#2{{\textstyle{#1 \over #2}}}
\def\lonep#1#2{{ ||#1||_1^{#2}}}
\def\dsp#1{{d_\star^{#1}}}
\def\dsmp#1{{(d_\star - 1)^{#1}}}

\def\half{{\textstyle{1\over2}}}
\def\quarter{{\textstyle{1\over4}}}

\def \kstarp#1 {{ K_{\star}^{#1}}}
\def \kst {{ K_\star^t}}
\def \Omega {{V}}
\def \rhon {{Z}}
\def \gtom {{\rightarrow_{m} \hspace{.015 in}  }}
\def \ks {{k^\star}}
\def \newmu {{\nu}}

\def \r {{\bf R}}
\def \var {{ \rm var }}
\def \r {{ \cal F}}

\def \la {{\langle}}
\def \ra {{\rangle}}
\def \b {{ \cal B}}
\def \kt {{ {\hat K } }}

\def \pmax {{\lambda}}

\def \shat {{\hat S}}

\def\r|{{\Bigr\vert}}
\def\l|{{\Bigl\vert}}

\def \R {{\bf R}}

\def\phi {\Phi}
\def\psistar {\psi_*}

\def\e{\epsilon}

\def \l {{W}}

\def\tmix{\tau_{\rm mix}}

\def\varepsilon{\mathchar"122 }

\def \one {{\mathbf 1}}
\def \chi {{\mathbf 1}}

\def\P {{ \bf P}}

\def\e {{ \bf {E}}}
\def\ds {{d_\star}}
\def\quartd{{  \alpha(\ds)}}

\def\Sq{{\cal S}_q}
\def\Sq-{{\cal S}_{q-1}}

\def\u{{\pi   }}

\def \lam {{\Lambda}}
\def \z {\bf Z}
\def\Px{\P_{\! \{x\}} \!}

\def\wS{\widetilde{S}}

\def\one{{\mathbf 1}}

\def\b{{\cal B}}

\def \et {{ \quartd}}
\def \lds {{ \cal L_{\star}}}
\newcommand{\be}{\begin{equation}}
\newcommand{\ee}{\end{equation}}
\newcommand{\lab}{\label}
\begin{document}
\title{The mixing time of the Thorp shuffle}
\author{
{\sc Ben Morris}\thanks{Department of Mathematics,
University of California, Davis.
Email:
{\tt morris@math.ucdavis.edu}.
This work was done while the author was at Indiana University and 
Microsoft Research.
} 
} 
\date{}
\maketitle 
\begin{abstract}
\noindent 
The Thorp shuffle is defined as follows. 
Cut the deck into two equal piles. 
Drop the first card from the left pile or the
right
pile according to the outcome of a fair coin flip;  
then drop from the other pile. 
Continue this way until both piles are empty.
We show that the mixing time for the  Thorp shuffle with $2^d$
cards is polynomial in $d$.
\end{abstract}
\setcounter{page}{1}
\section{Introduction} \lab{intro}
\subsection{The Thorp shuffle}
 How many shuffles are necessary to mix up a deck of cards? 
We refer to this as the {\it mixing time}
(see section \ref{mainresult} for a precise definition).
The mathematics of card shuffling has been studied extensively over the past 
several decades and most of the problems have been solved.  
Most famously, 
Bayer and Diaconis \cite{bd}
(in one of the few 
mathematical results to have made the front page of the New York Times) 
gave very precise bounds
for the Gilbert-Shannon-Reeds 
(riffle) shuffle model. 
Their bounds were correct even up to the constant factors.
For almost all natural shuffles matching upper and lower bounds are known 
(often even up to constants).
However, one card shuffling problem
has stood out for its resistance to attack. 

 In 1973, Thorp \cite{th} introduced the following shuffling procedure.
Assume that the number of cards, $n$, is even.
Cut the deck into two equal piles. 
Drop the first card from the left pile or the
right
pile according to the outcome of a fair coin flip;  
then drop from the other pile. 
Continue this way, 
with independent coin flips deciding whether to drop {\sc left-right} or 
{\sc right-left} each time, until 
both piles are empty. 

 The  Thorp shuffle, despite its simple description, 
has been hard
to analyze. The problem of determining its mixing 
time is, according to Persi Diaconis \cite{per},   
the ``longest-standing open card shuffling problem.''
It has long been  conjectured that
the mixing time is 
$O(\log^c n )$ for some constant $c$.
However, despite 
much effort
the only known upper bounds are
trivial ones of the form $O(n^c)$ that have circulated in the folklore.
The main contribution of this paper is to give the first 
poly log upper bound for the mixing time.

 We shall assume that the number of cards is
$2^d$ for a positive integer $d$. (Thus, our aim is to 
prove that the mixing time is polynomial in $d$.) In this case
the Thorp shuffle has a very appealing alternative description.
By writing the position of each card, from the bottom card ($0$)
to the top card ($2^{d} - 1$), in binary, we can view the 
cards as occupying the vertices of the 
$d$-dimensional unit hypercube $\{0,1\}^d$.
The Thorp shuffle proceeds in two stages. In the first stage,
an independent coin is flipped  for each edge $e$ in direction $1$
(i.e., each edge in the cube that connects two vertices that 
differ in only the first coordinate). 
If the coin lands heads, the cards at the endpoints of $e$ are 
interchanged; otherwise the cards remain in place. In the second stage, a 
``cyclic left bit shift'' is  performed for each card, where 
the card in position $(x_1, \dots, x_d)$ is moved to
$(x_2, \dots, x_d, x_1)$. 

We will actually use a 
slightly modified definition of the Thorp shuffle.
Say that an edge in the hypercube {\it rings} if its endpoints are
switched with probability $\half$. 
For $j = 1, \dots, d,$ let
$K_j$ be the transition kernel for the 
the process in which every edge $e$ in 
direction $j$ rings. 
\smallskip 
\\
\noindent {\bf Definition: Thorp shuffle.} 
The {\it Thorp shuffle} 
is the Markov chain whose transition kernel 
at time $n$ is $K_{j+1}$ if $j \equiv n \bmod d$.

Since $d$ iterations of this shuffle is equivalent to $d$
iterations of the shuffle described in \cite{th}, it is enough to
prove  a poly($d$) mixing
time bound for this new model.

 It is natural to consider the change in the deck after $d$ shuffles
 have been performed. (This 
represents one complete ``cycle''.) We will call this a {\it round}.
Using the language of network computing, 
a round of the Thorp shuffle is like passing the cards through $d$ levels
of a butterfly 
network (see, e.g., Knuth's book \cite{kb}),
where at each stage neighboring cards are 
interchanged with probability $\half$.
We note that in a recent breakthrough result,
\c{C}am \cite{hc} showed 
that the matrix $K_1 \cdots K_d K_1 \cdots K_{d-1}$
has strictly positive entries. 
This can be viewed as a result about the ``diameter'' of the 
Thorp shuffle; after a small number of steps there is a positive probability 
of being in any given state. 
However,
these probabilities are in general very small  so this does 
not imply a
good  bound for the mixing time.

 The main result of this paper is that indeed the mixing time is polynomial 
in $d$. Our proof uses 
{\it evolving sets}, a technique for bounding 
mixing times that was 
introduced by the author and Peres
in \cite{mp}.
Another paper that uses some of the same ideas 
is \cite{ex}, in which a variant of 
evolving sets is used to analyze the exclusion process. 
Evolving sets are related to the notion of strong stationary duality 
due to Diaconis and Fill \cite{df}.

\subsection{Statement of main result}
\lab{mainresult}
For a Markov chain on state space $V$ with uniform stationary 
distribution, 
define the {\it (uniform) mixing time}
by $$
 \tmix = 
 \min\Bigl\{n: \Bigl|  p^n(x,y)|V| -1 \Bigr| \leq \quarter
 \;\;\forall \,x,y \in \Omega               \Bigr\},$$
where $p^n(x,y)$ is the $n$-step transition probability from $x$ to $y$.
(This 
is a stricter definition of mixing time than the usual 
one involving total variation distance.)

Our main result is the following theorem.
\begin{theorem}
\lab{maintheorem}
The mixing time for the Thorp shuffle is
$O(d^{44})$. 
\end{theorem}

  In similar fashion to the analysis in \cite{mp}, we prove 
our mixing time bound based on an isoperimetric function 
we call the root profile. 
The paper is organized as follows. Following 
a brief introduction to evolving sets in section \ref{sec-es},
we devote much of the rest of the paper to proving a bound on the 
root profile.
In section \ref{l2tech} we show how 
$l^2$ techniques can be combined with evolving
sets to give a bound on the root profile. In section \ref{shuffling} we 
describe the chameleon process, a variant of 
the Thorp shuffle in which the cards have changing colors,
which is useful to bound mixing times.  
In section \ref{halfcards} we 
use the chameleon process to show 
that for a ``reversibilized'' version of the 
Thorp shuffle, any collection of cards 
(if viewed as indistinguishable)
mixes in
poly($d$) time. 
In section \ref{psisec}, we state the main 
technical result of this paper (proved in section
\ref{mlsec}), which says that the transition kernel
for the Thorp shuffle contracts functions in a certain 
$l^2$ sense; then we use this to
obtain our bound on the root profile. 
Next, armed with a good bound on the root profile 
we prove Theorem \ref{maintheorem} in section \ref{sec-mt}.
We conclude with proofs of some technical lemmas in sections 
and \ref{mlsec2} and \ref{mlsec}.
\section{Evolving sets}
\lab{sec-es}
\noindent
We  will now give a brief overview of 
evolving sets (see \cite{mp} for a more detailed account).
Let $\{p(x,y)\}$ be transition probabilities
for an irreducible, aperiodic Markov chain on a finite state space $V$.
Assume that the chain has a uniform stationary distribution 
(which means that $p$ is doubly stochastic:
$\sum_{x \in V} p(x,y)= 1$ for all $y \in V$).
For subsets $S \subset V$, 
define $p(S, y) := \sum_{x \in S} p(x,y)$. 
\smallskip 
\\
\noindent {\bf Definition: Evolving sets.} 
The {\it evolving set} process is 
the Markov chain $\{S_n\}$ on 
{\it subsets} of $V$ with the following transition
rule. If the current state $S_n$ 
is $S \subset V$, choose $U$ uniformly
from $[0,1]$ and let the next state $S_{n+1}$ be
$$
\wS=\{ y: p(S,y) \geq U \} \,. 
$$

 Write
$\P_S\Big( \cdot \Big):= \P\Big(\,\cdot\, \Bigl| \, S_0 = S \Big)$
and similarly for $\e_S\Big( \cdot \Big)$.
Evolving sets have the following properties (see \cite{mp}).
\begin{enumerate}
\item 
The sequence $\{|S_n|\}_{n \ge 0}$ forms a martingale. 
\item 
For all $n \geq 0$ and $x,y \in V$ we have
\[
p^n(x,y) = \Px\left(y \in S_n \right).
\]
\item
The sequence of complements
$\{S_n^c\}_{n \ge 0}$ 
is also an evolving set process, with the same transition probabilities.
\end{enumerate}

As in \cite{mp}, we will prove 
our mixing time bound using an isoperimetric quantity that we 
denote by $\psi$, which is defined as follows. 
For $S \subset V$, define
$$
\psi(S) := 1 - \e_S \sqrt{ \frac{ |\wS|}{ |S|}}.
$$
Define $\psi(x)$ for $x \in [0,1/2]$ by
\be \lab{defgam}
\psi(x) = \inf \{ \psi(S):  |S| \leq x |V| \},
\ee
and for $x>1/2$, let $\psi(x):=\psi_*=\psi(\half)$.
Observe that $\psi$ is non-negative and (weakly) decreasing on $[0,\infty)$.
We will call the function $\psi$ the {\it root profile.}

\section{From $\ell^2$ bounds to a bound on $\psi$}
\lab{l2tech}
In this section, we show how to use $l^2$ techniques to obtain a 
bound on the root profile.

Let $p(x,y)$ be a doubly stochastic
Markov chain on the state space $V$.
For functions $f: V \to [0,1],$ define 
$\lone{f} := {1 \over |V|} \sum_{x\in V} f(x)$ and
$\ltwo{f} := ({1 \over |V|} \sum_{x\in V} f(x)^2)^{1/2}$.
For $S \subset V$, define $\chi_S: V \to [0,1]$ by
\[
\chi_S(x) =
\left\{\begin{array}{ll}
1 &
\mbox{if $x \in S$;} \\
0 & \mbox{otherwise.}\\
\end{array}
\right.
\]

\begin{lemma} 
\label{secondcase}
Let $\wS$ be the next step in the evolving set process starting from $S$,
i.e., $\wS = \{y: p(S,y) > U\}$, where $U$ is uniform.
Let $\alpha = { \ss{ p(S, \, \cdot \,)} \over \lone{ \one_S}}$. 
Then
\[
\e\Bigl( \sqrt{ {|\wS| \over |S|}} \Bigr) \leq 
\Bigl[ { \alpha} (2 - 
\alpha)
\Bigr]^\quarter.
\]
\end{lemma}
\begin{proof}
Let $\Lambda$ be an independent copy of $\wS$, i.e., 
$\Lambda = \{y: p(S,y) > U'\}$, 
for an independent uniform random variable $U'$.
Note that either $\wS \subseteq \Lambda$ or 
$\Lambda \subseteq \wS$ (depending on which of the uniform variables 
$U,U'$ is larger). 
Let $X = |\wS \cap \lam|$ and 
$Y = |\wS \cup \lam|$. 
Then
\begin{eqnarray}
\label{etwo}
\Bigl[ \e( \sqrt{|\wS|}) \Bigr]^2 &=&
\e( \sqrt{|\wS||\Lambda|}) \\
&=& 
\e( \sqrt{XY}) \\
&\le& 
\sqrt{ \e(X) \e(Y)} \\
\lab{dec}
&=& \sqrt{ \e(X) (2 |S| - \e(X))},
\end{eqnarray}
where the first inequality is Cauchy Schwarz and the second inequality 
follows from the fact that $\e(X + Y) = 2 \e(\wS) = 2|S|$.
But
\begin{eqnarray}
\e(X) &=& \sum_{y \in V} \P( y \in \wS \cap \lam) \\
&=&  \sum_{y \in V} \P( y \in \wS)^2 \\
&=&  \sum_{y \in V} p(S, y)^2 = |V| \cdot \ss{ p(S, \cdot \,)},
\end{eqnarray}
so dividing the LHS of (\ref{etwo}) and the RHS of (\ref{dec}) by 
$|S| = |V| \cdot \lone{ \one_S}$ 
and then taking a square root yields the lemma. 
\end{proof}
\begin{remark}
\lab{foll}
The same proof shows that if $\wS = \{y: f(y) > U\}$ 
for $f:V \to [0,1]$ arbitrary, 
then
\[
\e\Bigl( \sqrt{ { \lone{\one_{\wS}} \over \lone{f}}} \Bigr) \leq 
\Bigl[ { \alpha} (2 - 
\alpha)
\Bigr]^\quarter,
\]
where 
$\alpha = 
{\ss{K^tf} \over \lone{f}}$ for $K$ the transition kernel.
Note also that if we define $\Delta := 1 - \alpha,$ then
\[
\Bigl[ { \alpha} (2 - 
\alpha)
\Bigr]^\quarter = (1 - \Delta^2)^\quarter \leq 1 - {\Delta^2 \over 4}.
\]
\end{remark}

\section{Chameleon process}
\lab{shuffling}

It will be convenient to study the card shuffle that behaves like the
Thorp shuffle for the first $d$ steps ($K_1, \dots, K_d$), and 
then like a ``reverse Thorp shuffle'' for the next $d$ steps
($K_d, \dots, K_1$). We will call this the zigzag shuffle. Every $2d$
steps of the zigzag shuffle will be called a round. (So a round 
of the zigzag shuffle is a round of the Thorp shuffle followed by a
round of a time-reversed Thorp shuffle.) 

Let $\alpha$ be large enough so that $4 \alpha^{-d} \leq 2^{-d-1}
4^{-d}$ for all $d \geq 1$ and let $c$ be an integer large enough so that
$[4e^{-c}]^d \beta \log \alpha \, cd^5
\le \alpha^{-d}$ for all $d \geq 1$, where 
$\beta = 2056 \cdot 64 \cdot 5$. 

The chameleon process is an extension of the zigzag shuffle.
The cards move in the same way as in the zigzag shuffle, but they 
also have colors, which can be red, white, black or pink. 
Initially, the
cards are colored as follows. 
There is a sequence of cards $x_1, \dots, x_b$ for some $b >
2^{d-1}$ such that cards $x_1, \dots, x_{b-1}$
are colored white, card $x_b$ is colored red, and
the 
remaining cards are colored black. The cards can change
color in two ways. The first way is called {\it pinkening}, which
takes place when an edge connecting a red card to a white
card rings; in this case both cards are re-colored pink. The second
way is called {\it de-pinking}, which takes place at the end of every
$64cd$
rounds of shuffling; 
in this case all of the pink cards are collectively
re-colored red or white, with probaility $\half$ each. (A process of
this type was first used in
\cite{ex} to analyze the exclusion process.) Note that black cards can 
never change color. 

Let $X_n$ be the zigzag shuffle. 
For $j = 1, \dots, 2^d$, we 
will write $X_n(j)$ for the
position of card $j$ at time $n$. 
If $S = \{z_1, \dots, z_k\}$ is a set of cards, define
$X_n(S) = \{X_n(z_1), \dots, X_n(z_k)\}$.
Let $W_n = X_n\Bigl( \{1, \dots, b \} \Bigr)$ be the {\it unordered} 
set of locations of nonblack (i.e., white, red or pink) cards at time $n$.
For vertices $x$ in the hypercube, define 
\[
\rho_n(x) = \one(\mbox{there is a red card at
  $x$ at time $n$}) + \half
\one(\mbox{there is a pink card at
  $x$ at time $n$}).
\]
The following lemma
indicates the fundamental relationship between the chameleon process
and the zigzag shuffle. 

\begin{lemma} 
\label{chamlemma}
Consider the chameleon process with $b$ nonblack
cards. Then
\[
\P\Bigl( X_n(x_b) = x \,\Bigl|\, W_1, W_2, \dots \Bigr)
= 
\e\Bigl( \rho_n(x) \,\Bigl|\, W_1, W_2, \dots \Bigr).
\]
\end{lemma}
\begin{proof}
We will use induction on $n$. The base case $n = 0$ is trivial because
there is initially only one red ball which is located at the position
of card $x_b$. Now assume that the result holds for $n$. Let $e$ be the
edge incident to $x$ that rings at time $n$ and let $x'$ be the
neighbor of $x$ across $e$. Let $A_1$, $A_2$ and $A_3$ be the events
corresponding to the following three possible values of 
$\Bigl(W_n \cap \{x,x'\}, W_{n+1} \cap \{x, x'\}\Bigr)$
when $x \in W_{n+1}$:
\begin{enumerate}
\item
$( \{x,x'\}, \{x, x'\})$;
\item
$(\{x'\}, \{x\})$; 
\item
$(\{x\}, \{x\})$.
\end{enumerate}
Let $\f_n = \sigma( \rho_n(x), \rho_n(x'))$.
Note that 
\be
\lab{coloreq}
\e\Bigl(\rho_{n+1}(x) \,\Bigl|\, \f_n, W_1, W_2, \dots \Bigr)
= \Bigl(\half \rho_n(x) + \half \rho_n(x')\Bigr) \one(A_1) +
\rho_n(x') \one(A_2) +
\rho_n(x) \one(A_3).
\ee
Define $\mu_n(\cdot) = 
\P\Bigl( X_n(x_b) = \,\cdot\, \,\Bigl|\, W_1, W_2, \dots \Bigr).$
Then
\be
\lab{probeq}
\mu_{n+1}(x) 
= \Bigl(\half \mu_n(x) + \half \mu_n(x')\Bigr) \one(A_1) +
\mu_n(x') \one(A_2) +
\mu_n(x) \one(A_3).
\ee
But by induction we have
\[ 
\mu_n(x) = \e\Bigl( \rho_n(x) \,\Bigl|\, W_1, W_2, \dots \Bigr);
\hspace{.3 in}
\mu_n(x') = \e\Bigl( \rho_n(x') \,\Bigl|\, W_1, W_2, \dots \Bigr).
\]
To complete the proof, take the conditional expectation given 
$W_1, W_2, \dots$ of both sides
of
(\ref{coloreq}) and  combine with equation (\ref{probeq}).
\end{proof}
\begin{remark}
\label{expectation1}
Note that 
\be
\lab{e1}
\e\Bigl( \sum_x \rho_n(x) \, \Bigl| \, W_1, W_2, \dots \Bigr)
= \sum_x \P\Bigl( X_n(b) = x \,\Bigl|\, W_1, W_2, \dots
\Bigr) = 1.
\ee
\end{remark}
\section{Indistinguishable cards mix in poly time}
\label{halfcards}
Let $\Lambda$ be a set of cards. Then
the process $\{X_n(\Lambda): n \geq 0\}$
is a Markov chain. The following lemma says that the uniform 
mixing time for this chain is $O(d^5)$.

\begin{lemma} There is a universal constant $b \in \z$ such that  
if $m = b d^5$ then
\label{pmaxlemma}
\be
\max_{\Lambda, \Lambda'} \,\, 
\Bigl|
{\textstyle {2^d \choose |\Lambda|}}
\P( \Lambda \to_m \Lambda') - 1 \Bigr|  \leq \quarter,
\ee
where we write $\Lambda  {\rightarrow_m} \Lambda'$ for the event that 
$X_m(\Lambda) = \Lambda'$.  
\end{lemma}
\begin{proof}
It is enough to consider sets $\Lambda$ with $|\Lambda| \geq 2^{d-1}$.
(Otherwise, consider $\Lambda^c$.)
Let $\alpha, \beta$ and $c$ be defined as in section \ref{shuffling},
let $b \ge \beta \log \alpha \, c,$ and let 
$m = bd^5$. 
For $j \in \{1, \dots, 2^d\}$ define 
\[
\pmax(j) = 
\max_{|S| = j}
\max_{|S'| = j}
 \Bigl|   {\textstyle{2^d \choose j}}
\P(S \gtom S')  - 1
\Bigr|.
\]
We will show that
for all $k \geq 2^{d-1}$, we have 
\be
\lab{pme}
\pmax(k) \leq 
\ks 4^{-d} ,
\ee
where $\ks = 2^d - k$. 
This yields the lemma because the r.h.s.~of
(\ref{pme}) is at most $\quarter$ for all $d \geq 1$. 

Let $A$ and $B$ be disjoint sets of cards.
For $x \in A$, say that $x$ is {\it antisocial 
in round $j$} of the zigzag shuffle if 
at no point in round $j$
does an edge connecting $x$ to
a card in $B$ ring. Let $Z(A,B,j)$ denote the number of cards that 
are antisocial in round $j$. 
We say that
$A$ {\it avoids} $B$ if $Z(A,B,j) > \sfrac{7}{8} |A|$ for 
$64cd$ consecutive rounds $j$ before time $m$. If $S$ is a set of cards, say
that {\it $S$ mixes}
if
there do not exist disjoint sets $A, B$ of cards with $|A| \leq \half |S|$
and $A \cup B = S$ such that $A$ avoids $B$.

We will verify (\ref{pme}) 
by induction on $\ks$.
The base case $\ks = 0$ ($k = 2^d$) is trivial. Suppose it's true for $k$, 
where $k > 2^{d-1}$ and
consider $k-1$. 
Fix a set of cards $S = \{x_1, \dots, x_k\}$ and consider the 
corresponding chameleon process.
Let $\f = \sigma(X_n(S) : n \geq 0)$.
Let  $\rhon_n = \sum_x \rho_{64cd^2n}(x)$ be the
total amount of ``red paint'' in the system after $64cdn$ rounds 
of the chameleon process. 
Define $\rhon_n^\sharp = \min( \rhon_n, k - \rhon_n)$. 
Note that $\lim_{n \to \infty} \rhon_n^\sharp = 0$ a.s.

 Fix $n$ such that $64cd^2 n \leq m$, and
let $A_n$ be either the set of cards that are red 
or the set of cards that are white at the start of round
$64cdn$, 
according to whether $\rhon_n \leq k/2$ or $\rhon_n > k/2$, 
respectively. 
Let $P$ denote the number of cards pinkened during
the next $64cd$ rounds.
Let $B_n = S - A_n$.
When $S$ mixes, 
$A_n$ doesn't avoid $B_n$.
We claim that this ensures that $P \geq {|A_n| \over 8d}$.
Consider a round $j$ such that $Z(A_n, B_n, j) \geq 
\sfrac{7}{8} |A_n|$.
Note that after an edge connecting 
a card $x$ in $A_n$ to a card $y$ in $B_n$ rings, 
at least one of the resulting
cards is pink. Let us
associate that pink card with $x$. 
(If both endpoints are pink then choose one of them arbitrarily.)
Since at least a fraction $1/8$ 
of the cards in $A_n$ will have a pink card associated to them in this
round, and since any given pink card can be associated to at most $d$ 
cards in $A_n$ in this round, the number of pink cards at the end of
this round must be at least ${|A_n| \over 8d}$. 
It follows that $P \geq {|A_n| \over 8d}$.

 Note that $\rhon_{n+1}$ is either $\rhon_n + \half P$ or 
$\rhon_n - \half P$, with probability $\half$ each. Thus,
if we write $E$ for the event that $S$ does not mix, then
\begin{eqnarray}
\e \Bigl( \sqrt{ \rhon_{n+1}^\sharp}
\, \Bigl| \, P, \rhon_n, \f, E^c \Bigr)
 &=&
\e\Bigl(
\half \sqrt{ (\rhon_n + \half P)^\sharp} + 
\half \sqrt{ (\rhon_n - \half P)^\sharp}
\, \Bigl| \, Z_n, \f, E^c \Bigr)
 \\ 
&\leq&
\sqrt{ \rhon_n^\sharp } \, { \sqrt{ 1 + \sfrac{1}{16d}} + \sqrt{1 - 
\sfrac{1}{16d}} \over
  2} \\
&\leq& \sqrt{ \rhon_n^\sharp } \, \exp\Bigl[{-{1 \over  2056d^2}}\Bigr],
\end{eqnarray}
where the first inequality follows from the concavity of the square
root, 
and the second inequality follows from the fact that
$\half \sqrt{1 + u} +  \half \sqrt{1 - u} \leq \exp( - u^2/8)$
whenever $u \in [0,1]$ (see \cite{mp}, Lemma 9).

Thus, since $\rhon_0 = 1$, it follows that
\be
\lab{expdecay}
\e \Bigl( \sqrt{ \rhon_{n}^\sharp}
\, \Bigl| \, \f, \mbox{$S$ mixes} \Bigr) \leq 
\exp\Bigl[{-{n \over  2056d^2}}\Bigr]
\ee
for all $n$. Define $Z_\infty = \lim_{n \to \infty} Z_n$. 
(Note that for any $S'$ we have 
$\e(Z_\infty \,|\, S \gtom S') = 1$; see the remark immediately
following 
Lemma \ref{chamlemma}.)
Lemma \ref{chamlemma} 
implies that for all $y \in S'$ we have
\begin{eqnarray}
\Bigl| \P\Bigl( 
X_m(x_k) = y \, \Bigl| \, S \gtom S'\Bigr)
- \sfrac{1}{k} \Bigr|
&=&
\Bigl|
\e \Bigl( \rho_m(y) 
- \sfrac{1}{k} Z_\infty 
\, \Bigl| \, S \gtom S' \Bigr)
\Bigr|
 \\
&\leq& \e \Bigl(
| \rho_m(y) - \sfrac{1}{k} Z_\infty | \,\Bigl|\,
S \gtom S'
\Bigr)
 \\
\label{combine1}
&\leq& 
\P(\rho_m \notin \{0,k\} 
\, | \, S \gtom S').
\end{eqnarray}
Let $E$ be the event that $S$ does not mix. 
Lemma \ref{edef} in Appendix A gives $\P(E \,|\, S \gtom S') 
\leq \alpha^{-d} {1 + \pmax(k) \over 1 - \pmax(k)}$. 
Hence
\begin{eqnarray}
\P\Bigl(\rho_m \notin \{0,k\} \, \Bigl| \, S \gtom S'\Bigr) &\leq&
\P\Bigl(E \,\Bigl|\, S \gtom S'\Bigr) + 
\P\Bigl(\rho_m \notin \{0,k\} \, \Bigl| \, S \gtom S', E^c\Bigr) \\
&\leq&
\alpha^{-d} {1 + \pmax(k) \over 1 - \pmax(k)} + 
\P\Bigl(\rho_m \notin \{0,k\} \, \Bigl| \, S \gtom S', E^c\Bigr) \\
&\leq&
\label{combine2}
3 \alpha^{-d} +
\P\Bigl(\rho_m \notin \{0,k\} \, \Bigl| \, S \gtom S', E^c\Bigr), 
\end{eqnarray}
where the third inequality holds because $\pmax(k) \leq \quarter$
by induction. But
\begin{eqnarray}
\P\Bigl(\rho_m \notin \{0,k\} \, \Bigl| \, S \gtom S', E^c\Bigr)
&\leq&
\e\Bigl(\rhon_{m/64cd^2}^\sharp \, \Bigl| \, S \gtom S', E^c \Bigr) \\
&\leq& \exp\Bigl[{-{m \over  2056 \cdot 64 \cdot c d^4}}\Bigr]
\label{combine3}
\leq \alpha^{-d},
\end{eqnarray}
where the second inequality follows from equation (\ref{expdecay}). 
Combining equations 
(\ref{combine1}), 
(\ref{combine2}), and 
(\ref{combine3}) gives 
\begin{eqnarray}
\label{onecard}
\Bigl| 
\P\Bigl( 
X_m(x_k) = y \, \Bigl| \, S \gtom S'\Bigr)
-
\sfrac{1}{k} \Bigr| &\leq& 
4 \alpha^{-d}.
\end{eqnarray}

Now fix a set of cards $\Lambda$ with $|\Lambda| = k-1$ and 
 and let $z \notin \Lambda$.
Define $\Lambda_z = \Lambda \cup z$.
Fix a set $\Lambda'$ of vertices of the hypercube
with $|\Lambda'| =  k-1$. For $w \notin \Lambda'$, define
\begin{eqnarray}
x_w = \P\Bigl( \Lambda_z \gtom \Lambda'_w \Bigr)
\hspace{.3 in}
\Delta x_w = x_w - \choosek \\
y_w = \P\Bigl( z \gtom w \, \Bigl| \, \Lambda_z \gtom \Lambda'_w \Bigr)
\hspace{.45 in}
\Delta y_w = y_w - 1/k.
\end{eqnarray}
Note that $\Bigl| \{w: w \notin \Lambda'\} \Bigr| = \ks + 1$,
and ${\ks + 1 \over k} {\textstyle {2^d \choose k}^{-1}}
= 
{\textstyle {2^d \choose k-1}^{-1}}$. It follows that
\begin{eqnarray}
\Bigl| \P(\Lambda \gtom \Lambda') - 
{\textstyle {2^d \choose k-1}^{-1}} 
\Bigr|
&=&
\Bigl|
\sum_{w \notin \Lambda'} 
\P\Bigl( \Lambda_z \gtom \Lambda'_{w}, z \gtom w \Bigr) 
- {\textstyle{1 \over k}}{\textstyle {2^d \choose k}^{-1}}
\Bigr|
\\
&=& 
\Bigl|
\sum_{w \notin \Lambda'} 
x_w y_w
- {\textstyle{1 \over k}}{\textstyle {2^d \choose k}^{-1}}
\Bigr|
\\
\label{pretri}
&=&
\Bigl |
\sum_{w \notin \Lambda'}
\Delta x_w {\textstyle {1 \over k}}
+
\Delta y_w \choosek
+
\Delta x_w 
\Delta y_w 
\Bigr | .
\end{eqnarray}
Note that
\be
\lab{firstin}
|\Delta x_w| \leq \ks 4^{-d} \choosek
\leq \choosek,
\ee
where the first inequality is induction and 
the second inequality holds because $\ks \leq 2^d$.
Also, 
equation (\ref{onecard}) implies that
\be
\lab{twoin}
| \Delta y_w | \leq 4 \alpha^{-d} 
\leq {\textstyle {1 \over 2k}}
4^{-d},
\ee
for all $d \geq 1$ by the definition of $\alpha$. 
Thus, using equations (\ref{firstin}),(\ref{twoin}) 
and the triangle inequality,
equation (\ref{pretri})
becomes
\begin{eqnarray*}
\Bigl| \P(\Lambda \gtom \Lambda') - 
{\textstyle {2^d \choose k-1}^{-1}} 
\Bigr|
&\le& 
{\ks + 1 \over k} \Bigl[
 \ks 4^{-d} \choosek +  4^{-d} \choosek
\Bigr] \\
&=&
{1 \over k} 
(\ks + 1)^2 
4^{-d}
\choosek \\
&=& (\ks + 1) 4^{-d} \choosekm = (k - 1)^\star 
4^{-d}
\choosekm.
\end{eqnarray*}
Since this is true for all $\Lambda$ with $|\Lambda| = k-1$
the proof is complete. 
\end{proof}
Let $K$ be the transition kernel for one round of the Thorp shuffle,
and 
let $K^t$ be the transpose of $K$, defined by
$K^t(x,y) = K(y,x)$. Note that $K^t$ is the time-reversal of $K$. 
Let $\kt := KK^t$ be the transition kernel for one round of the zigzag
shuffle. Let $\{Z_n : n\geq 0 \}$ be a Markov chain with transition
kernel $\kt$. 
Then
Lemma \ref{halfcards} implies that for any set of cards $B$,
the uniform mixing time for the process $\{Z_n(B) : n \geq 0\}$ is at most
$bd^4$. Thus, using standard facts about geometric convergence
and the uniform mixing time,
we can conclude that for a universal constant $C$ we
have 
\be
\lab{fred}
\max_{B'} {\textstyle {2^d \choose |B|}} \P\Bigl( Z_{kCd^4}(B) = B'\Bigr) 
\leq 1 + e^{-k},
\ee
for all $k \geq 1$. \\
\\
{\bf Truncated Thorp shuffle.}
Fix $\ds \leq d$. Define the {\em $\ds$--truncated Thorp shuffle}
as the Markov chain with transition kernel
$K_\star = K_1 \dots K_\ds$. This is a ``partial round'' of the Thorp shuffle,
with  steps
$\ds+1$ through $d$  
censored. 
To make things irreducible, we define the state space as the set of
states reachable from an (arbitrary) fixed starting state. 

 Define the $\ds$--truncated zigzag shuffle as the Markov chain
with transition kernel $K_\star \kst$.
Note that
we can think of this shuffle as 
a product of $2^{d - \ds}$ copies of a ``$\ds$-dimensional''
zigzag shuffle, where the cards occupy $2^{d - \ds}$ 
(disconnected) hypercubes of dimension
$\ds$. 
Combining this observation with equation (\ref{fred})
yields the following corollary to Lemma \ref{halfcards}.

\begin{corollary}
\label{trunccor}
Fix $\ds \geq 2$ and let $\{Z_n : n\geq 0\}$ be the $\ds$--truncated 
zigzag shuffle. 
There is a universal constant $c$ such that
if $l = kcd (\ds - 1)^4$, then
\be
\max_{B'} {\textstyle {2^d \choose |B|}} \P\Bigl( Z_l(B) = B'\Bigr) 
\leq \exp(\exp(-k)),
\ee
for all $k \geq 1$. 
\end{corollary}
\begin{proof}
Let $c = 2^5C$. Then $l \geq 2kdC\dsp{4}$, so
equation (\ref{fred}) implies that
\begin{eqnarray*}
\max_{B'} {\textstyle {2^d \choose |B|}} \P\Bigl( Z_{l}(B) = B'\Bigr) 
&\leq& (1 + e^{-2k d})^{2^{d- \ds}} \\
&\leq& \exp(2^d\exp(-2dk)) \\
&\leq& \exp(\exp(-k)),
\end{eqnarray*}
for all $d \geq 1$. 
\end{proof}
\section{A bound on the root profile}
\lab{psisec}
We will need the following technical result, which is proved in 
Appendix B.\\
\\
{\noindent 
{\bf Corollary \ref{cl}} 
{\em 
Fix $S \subset \Omega$ and let $
x = {|S| \over (2^d)!} = \lone{\one_S}$.
Let $p(\,\cdot\,, \, \cdot \,)$ be the transition kernel for one round 
of the Thorp shuffle.
Then there is a universal constant $C> 0$ such that 
\[
\ss{p(S, \cdot)} \leq x^{1 + {C/d^{14} }}.
\]
}}
We are now ready to obtain a bound on 
the root profile of the Thorp shuffle. 
\begin{lemma}
\lab{psi}
Let $\psi$ be the root profile of the Markov chain which each 
step performs a round of the Thorp shuffle ($K_1 K_2 \cdots K_d$).
There is 
a universal constant $c > 0$ such that
\be
\lab{maxe2}
\psi(x) \geq \max\Bigl( 1 - x^{c / 2d^{42}}, {c d^{-28}}\Bigr).
\ee
\end{lemma}
\begin{proof}
Let $C$ be the constant appearing in Corollary \ref{cl}.
We will show that there is a universal constant $B > 0$ such that
\be
\lab{maxe}
\psi(x) \geq \max\Bigl( 1 - x^{CB/2d^{42}}, {B d^{-28}}\Bigr).
\ee
Setting $c = \min(BC, C)$ will then yield the lemma. 
First, we show that $\psistar \geq {B d^{-28}}$.
Fix $S$ with ${ |S| \over (2^d)!} = x \leq \half$ and let
\[
\wS = \{y: p(S, y) > U\},
\]
where $\{p(x,y)\}$ are the transition probabilities for one round of
the Thorp shuffle.
The remark following Lemma 
\ref{secondcase} implies that
\[
\e \sqrt{|\wS^\sharp| \over |S^\sharp|}
\leq 1 - {\Delta^2 \over 4},
\] 
where $\Delta = 1 - { \ss{p(S, \,\cdot\,)} \over \lone{\one_S}}$, and 
Corollary \ref{cl} implies that
$\ss{ p(S, \, \cdot \,)} \leq x^{C/d^{14}} \lone{\one_S} \leq
2^{-C/d^{14}} \lone{\one_S}$. Thus
\begin{eqnarray}
\Delta &\geq& 1 - 2^{- {C  d^{-14}}} \\
&=& 1 - e^{- {C \log 2 \,  d^{-14}}} \\
&\geq& {A d^{-14}},
\end{eqnarray}
for a universal constant $A>0$, and hence $1 - {\Delta^2 \over 4} \leq 
1 - {B d^{-28}}$ for a universal constant $B \in (0,\quarter)$.
(The fact that we can take $B < \quarter$ will be used later on.) 
Since this holds for all $S$ with $|S| \leq \half (2^d)!$, we conclude that
$\psistar \geq {B d^{-28}}$. To complete the proof of Lemma \ref{psi}, 
we must 
show that equation (\ref{maxe}) holds when
the max is achieved by the first term. Suppose that 
$1 - x^{CB/2d^{42}} \geq {B d^{-28}}$. Then
\be
\lab{other}
x \leq (1 - Bd^{-28})^{2d^{42}/CB} \leq \exp(- 2C^{-1} d^{14}).  
\ee
Assume that (\ref{other}) holds. Lemma \ref{secondcase} gives
\be
\lab{from}
\e \sqrt{|\wS^\sharp| \over |S^\sharp|}
\leq (\alpha (2 - \alpha))^\quarter
\leq (2 \alpha)^\quarter,
\ee
where $\alpha = 
{\ss{ p(S, \, \cdot \,)} \over \lone{\one_S}}$.
Equation (\ref{other}) implies that
\[
x^{C/2d^{14}} \leq e^{-1} < {\half},
\]
and hence
\be
\lab{twobound}
2 \leq x^{- {C/2d^{14}}}.
\ee
Furthermore, Corollary \ref{cl} implies that
$\alpha 
 \leq x^{C/d^{14}}$.
Plugging this and (\ref{twobound}) into
(\ref{from}) gives
\be
\lab{to}
\e \sqrt{|\wS^\sharp| \over |S^\sharp|}
\leq (x^{{-C/2d^{14}}} x^{C /d^{14}}    )^\quarter
= x^{C/ 8d^{14}}
\le x^{CB /2d^{42}},
\ee
since $B < \quarter$ (and $x \leq 1$). 
\end{proof}

\section{Proof of main result}
\lab{sec-mt}
\begin{proofof}{Proof of Theorem \ref{maintheorem}}
We shall start by bounding the mixing time of the Markov chain that 
does an entire round of the Thorp shuffle each step.
Recall that 
the root profile $\psi: [0, \infty) \to \R$ is 
defined by
\[
\psi(x) =
\left\{\begin{array}{ll}
\inf \{ \psi(S): |S| \leq x |\Omega| \} &
\mbox{if $x \in [0, \half]$;} \\
\psi_*    & \mbox{if $x > \half$,}\\
\end{array}
\right.
\]
where $\psi_* = \psi(\half)$.
Thus  $\psi$ is (weakly) decreasing on $[0,\infty)$.

Let $h(z) := 1 - \psi(1/z^2)$. 
Since $\psi(x) = \psistar$ for all real 
numbers $x \geq \half$, the function $h$ is 
well-defined even for $z \leq 1$.
Note that $h$
 is nonincreasing. 
In \cite{mp} it is shown (see section 5 and the part of section
3 entitled ``Derivation of Theorem 1 from Lemma 3 and Theorem 4'') that 
there
is a 
sequence of random variables $\{Z_n : n \geq 0\}$
that satisfies 
$Z_0 = \sqrt{ |\Omega|}$ and 
\begin{eqnarray} \lab{pref}
\e\left( \frac{Z_{n+1}}{Z_n} \Big| Z_n\right) 
\leq h(Z_n), \label{fzero}
\end{eqnarray}
such that
\be
\lab{mixing}
\tmix \le 2 \min\{n: \e(Z_n) \le \half \}.
\ee
Lemma \ref{psi} gave the following 
bound on the root profile: 
\be
\lab{maxe2}
\psi(x) \geq \max\Bigl( 1 - x^{c / 2d^{42}}, {c d^{-28}}\Bigr),
\ee
for a universal constant $c > 0$.
Thus $h \leq g$, where $g$ is defined by
\[
g(z) = \min\Bigl(z^{-c / d^{42}}, 1 - {c d^{-28}}\Bigr), 
\]
and hence $\e( Z_{n+1} | Z_n) \le g(Z_n) Z_n$. 
Let $f(z) = z g(z) = 
\min\Bigl(z^{1-{c / d^{42}}}, z(1 - {c d^{-28}})\Bigr)$. 
Note that $f$ is increasing and, as the minimum of two concave functions,
is concave.
We claim that $\e(Z_n) \leq f^n(Z_0)$, where $f^n$ is the $n$-fold 
iterate of $f$. We verify this by induction. The base case 
$n = 0$ is immediate. Suppose that the claim holds for $n$. Then
\begin{eqnarray}
\e(Z_{n+1}) &=& \e ( \e( Z_{n+1} | Z_n )) \\
&\le&  \e( f(Z_n)) \\
&\le& f( \e(Z_n)) \\
&\le& f( f^n(Z_0)) = f^{n+1}(Z_0),
\end{eqnarray} 
where the third line follows from concavity and the last line 
is the induction hypothesis.
Let 
\[
f_1 =
z^{1-{c / d^{42}}}; \hspace{.3 in}
f_2 = z(1 - {c d^{-28}}),
\]
so 
that $f = \min(f_1, f_2)$. Then for all $m,n$
we have
\[
\e(Z_{m+n}) \leq f^{m+n}( \sqrt{Z_0}) \leq f_2^m( f_1^n( Z_0)).
\]
But $f_1^n(z) = z^{(1 - {c / d^{42}})^n} \le 
z^{ \exp( { -cn / d^{42}})}$, and 
$Z_0 = \sqrt{ |\Omega|} \leq (2^d)^{2^d} = 2^{d2^d}$. 
Thus, choosing $n \geq c^{-1}d^{43}$ gives
\[
f_1^n(Z_0) \le 2^{d 2^d e^{-d}},
\]
which is at most $4$ for all $d \geq 1$. 
Finally, since
\[
f_2^m(z) = z \Bigl(1 - {c d^{-28}}\Bigr)^m \le z e^{- {cm / d^{28}}}, 
\]
we have $f_2^m(4) \le 4 e^{- {cm/d^{28}}},$ which is at most 
$\half$ whenever $m \geq {c^{-1}d^{28} \log 8}.$ Putting this together,
we conclude 
that $\tmix \leq {2c^{-1}} ( d^{43} + d^{28} \log 8) = O(d^{43})$.
Since each round corresponds to $d$ Thorp shuffles we conclude that
the mixing time for the original model is $O(d^{44})$. 
\end{proofof}  

\section{Appendix A}
\label{mlsec2}
In this section we prove some large deviation results needed
in section \ref{halfcards}. We will adopt the notation of that 
section; 
for the convenience of the reader, we now give a brief recap. 
Let $A$ and $B$ be disjoint sets of cards.
For $x \in A$, say that $x$ is {\it antisocial 
in round $j$} of the zigzag shuffle if 
at no point in round $j$
does an edge connecting $x$ to
a card in $B$ ring. Let $Z(A,B,j)$ denote the number of cards that 
are antisocial in round $j$. 
We say that
$A$ {\it avoids} $B$ if $Z(A,B,j) > \sfrac{7}{8} |A|$ for 
$64cd$ consecutive rounds $j$ before time
$m$. If $S$ is a set of cards, say
that {\it $S$ mixes}
if
there do not exist disjoint sets $A, B$ of cards with $|A| \leq \half |S|$
and $A \cup B = S$ such that $A$ avoids $B$.

\begin{lemma}
\label{momentlemma}
Let $\{X_n : n \geq 0\}$ be the 
zigzag shuffle. 
Let $Z = Z(A, B, 1)$ be the number of cards that are
antisocial in the first round.
Define $\f_B = \sigma( X_1(B), 
\dots, X_d(B))$. 
Let 
$p = 1 - {|B| \over 2^d}$ 
and let $k = |A|$. 
For $\theta \geq 0$ 
define $\phi_p(\theta) = 1-p + p e^\theta$. 
Then for all $\theta \geq 0$ we have
\be
\lab{momineq}
\e \Bigl( e^{ \theta Z} \, \Bigl| \, \f_B \Bigr) \leq \phi_p(\theta)^k.
\ee
\end{lemma}
\begin{proof}
We verify this by induction on $d$. If $d = 1$ then the LHS of 
(\ref{momineq}) is $1$ if $p < 1$, and $e^{\theta k}$
otherwise, so (\ref{momineq}) holds. Now suppose that $d > 1$. 
Let $A'$ be the set of cards in $A$ not adjacent to $B$ in 
direction $1$, and let $k' = |A'|$.
Let $l$ be half the number of cards
in $A'$ adjacent to another card in $A'$ in direction $1$.
(Note that $l$ is an integer.) 
Let $k_0$ and $k_1$ be the
number of cards in $A'$ that end up 
with a leading $0$ and $1$, respectively, 
after the first step of the round
(i.e., after the edges in direction $1$ ring). Of those in the first
group, let $Z_0$ be the number that are antisocial, with a similar
definition for $Z_1$. 
Note that given $\f_B$, the random variables 
$k_0$ and $k_1$ are both distributed like
$W + l$, where $W \sim$ Binomial($k' - 2l$, $\half)$, and note 
that $Z = Z_0 + Z_1$. 
By induction, we have
\begin{eqnarray*}
\e\Bigl( e^{\theta Z} \, \Bigl| \, \f_B, X_1(A) \Bigr)
&=& 
\e\Bigl( e^{\theta Z_0} \, \Bigl| \, \f_B, X_1(A) \Bigr)
\e\Bigl( e^{\theta Z_1} \, \Bigl| \, \f_B, X_1(A) \Bigr) \\
&\leq& \phi_{p_0}(\theta)^{k_0}
\phi_{p_1}(\theta)^{k_1},
\end{eqnarray*}
where $p_0$ is the fraction of locations of 
the part of the hypercube with a leading $0$ not occupied by 
a card in $B$ after the first step, with a similar definition for
$p_1$.  It follows that 
$\e\Bigl( e^{\theta Z} \, \Bigl| \, \f_B, k_0, k_1 \Bigr)
\leq
\phi_{p_0}(\theta)^{k_0}
\phi_{p_1}(\theta)^{k_1}$. 
Hence
\begin{eqnarray*}
\e\Bigl( e^{\theta Z} \, \Bigl| \, \f_B \Bigr)
&\leq&
\sum_{i=0}^{k' - 2l} (\half)^{k' - 2l}
{k' - 2l \choose i}
\phi_{p_0}^{i}(\theta)
\phi_{p_1}^{k' - 2l - i}(\theta) 
\phi_{p_0}^{l}(\theta)
\phi_{p_1}^{l}(\theta)   \\
&=&
\Bigl[ 
\half \phi_{p_0}(\theta) + 
\half \phi_{p_1}(\theta)  \Bigr]^{k' - 2l} 
\phi_{p_0}^l(\theta)
\phi_{p_1}^l(\theta) \\
&\leq&
\Bigl[ 
\half \phi_{p_0}(\theta) + 
\half \phi_{p_1}(\theta)  \Bigr]^{k'}
= \phi_p(\theta)^{k'}, 
\end{eqnarray*}
where the last inequality follows from 
the AM-GM inequality and the final equality holds because
$p = \half (p_0 + p_1)$. This yields the lemma because $k' \leq k$.
\end{proof}
Lemma \ref{momentlemma} easily gives the following large deviation 
inequality. 
\begin{corollary}
\label{corused}
Suppose that $p \leq 3/4$. 
Then  
\[
\P\Bigl(Z > \sfrac{7}{8} k \, \Bigl| \, \f_B \Bigr) 
< e^{-k/64}.
\] 
\end{corollary}
\begin{proof}
We have
\begin{eqnarray}
\e( e^{ \theta(Z - pk)} \,|\, \f_B) &=&
e^{-pk \theta} \e( e^{\theta Z} \,|\, \f_B) \\
\label{mom}
&\leq& \Bigl[ (1 - p) e^{-p \theta} + p e^{\theta(1-p)} \Bigr]^k,
\end{eqnarray}
by Lemma \ref{momentlemma}.
The quantity inside the square brackets is $\e\Bigl( e^{\theta (Y-p)}
\Bigr)$, for a Bernoulli($p$) random variable $Y$. The inequality
$\e(e^{W}) \leq e^{\var(W)}$, valid when $\e(W) = 0$ and $W \leq 1$
(see, e.g., \cite{st}), 
implies that the quantity (\ref{mom}) is at most
$\exp\Bigl(\quarter \theta^2k\Bigr)$ 
if $\theta \le 1$. Letting $\theta = \quarter$
gives
\be
\lab{ebound}
\e\Bigl( \exp [\quarter(Z - pk)] \Bigr) \leq 
e^{k/64},
\ee
and hence
\begin{eqnarray}
\P\Bigl( Z > \sfrac{7}{8} k \,\Bigl|\, \f_B \Bigr) &=&
\P\Bigl( \exp [\quarter(Z - pk)] >
\exp\Bigl[ \sfrac{7k}{32} -
\sfrac{pk}{4}
\Bigr] \, | \, \f_B
\Bigr) \\
\lab{twopos}
&\leq&  
\exp\Bigl[ - \sfrac{7k}{32} +
\sfrac{pk}{4}
\Bigr]
\exp \Bigl[\sfrac{k}{64} \Bigr],
\end{eqnarray}
by Markov's inequality. Finally, since $p \leq 3/4$, the quantity 
(\ref{twopos}) is at most $e^{-k/64}$.
\end{proof}

The following lemma was used in the proof of Lemma
\ref{pmaxlemma} in section \ref{halfcards}.

\begin{lemma}
\lab{edef}
Fix a set of cards $S$ with $|S| \geq 2^{d-1}$. 
Then for any set $S'$ of vertices of the hypercube
we have
\[
\P\Bigl(\mbox{$S$ does not mix} \,\Bigl|\, S \gtom S'\Bigr) \leq \alpha^{-d}
{1 + \pmax(|S|) \over 1 - \pmax(|S|)}.
\]
\end{lemma}
\begin{proof}
Let $E$ be the event that $S$ does not mix.
We have
\begin{eqnarray*}
\P\Bigl(E, S \gtom S'\Bigr) &\leq&
\sum_{k \le \half |S|} \,\, 
\sum_{A: |A| = k}
\P\Bigl(\mbox{$A$ avoids $B$}, S \gtom S' \Bigr) \\
&\leq& 
2^{d-1} \max_k \Bigl[ 2^{dk} \max_{A: |A| = k} 
\P\Bigl(\mbox{$A$ avoids $B$}, S \gtom  S' \Bigr)\Bigr],
\end{eqnarray*}
where in the summations we write $B$ for $S - A$, 
the $2^{d-1}$ is an upper bound on the number of 
$k \leq \half |S|$, and the $2^{dk}$ is an upper bound on 
the number of sets $A$ with $|A| = k$.
Since $|A| \leq \half |S|$ and 
$A \cup B = S$, we must have $|B| \geq \quarter 2^{d}$.
Hence if 
$|A| = k$ then
\begin{eqnarray}
\P\Bigl( \mbox{$A$ avoids $B$}, S \gtom S' \Bigr) 
&\leq&
\sum_{B' \subset S'}
\P\Bigl( \mbox{$A$ avoids $B$}, B \gtom B', A \gtom A' \Bigr) \\
&\leq& 
\label{bon}
\sum_{B' \subset S'}
\P(B \gtom B') 
\P\Bigl( \mbox{$A$ avoids $B$} \, \Bigl| \, B \gtom B' \Bigr),
\end{eqnarray}
where in the summations, we write $A'$ for $S' - B'$.
But 
\be
\lab{sandy}
\P\Bigl( \mbox{$A$ avoids $B$} \, \Bigl| \, B \gtom B' \Bigr)
\leq \sum_{i=0}^{m}
\prod_{j=i}^{i+64cd-1} \P\Bigl(Z(A,B,j) > \sfrac{7k}{8}
\, \Bigl| \, B \gtom B'
\Bigr)
\leq m \Bigl( e^{-k/64} \Bigr)^{64cd},
\ee
where the last inequality follows from Corollary \ref{corused}.
Hence
\begin{eqnarray*}
\P\Bigl(\mbox{$A$ avoids $B$}, S \gtom  S' \Bigr)
&\leq&
\sum_{B' \subset S'}
{\P(B \gtom B')}
m e^{- ckd} \\
&\leq& 
2^{dk} m e^{- ckd} \max_{B'} {\P(B \gtom B')},
\end{eqnarray*}
where the $2^{dk}$ is an upper bound on the number of subsets 
$B' \subset S'$. 
But for any $B'$ we have
\[
\P(B \gtom B') \le \sum_{\shat: \shat \supset B'}
\P(S \gtom \shat) 
\le 
2^{dk} {\textstyle {2^d \choose |S|}^{-1}} 
(1 + \pmax(|S|)).
\]
It follows that 
\begin{eqnarray}
\P\Bigl(\mbox{$A$ avoids $B$}, S \gtom  S' \Bigr)
&\leq&
4^{dk}
m e^{-cdk} {\textstyle {2^d \choose |S|}^{-1}} 
(1 + \pmax(|S|)) \\
&=& 
[4e^{-c}]^d \beta \log \alpha \, c d^5
{\textstyle {2^d \choose |S|}^{-1}} 
(1 + \pmax(|S|)) \\
\label{later}
&\le& \alpha^{-d} {\textstyle {2^d \choose |S|}^{-1}} 
(1 + \pmax(|S|)), 
\end{eqnarray}
where the second inequality follows from the definition of $c$. 
Finally, since $\P(S \gtom S') \geq 
{\textstyle {2^d \choose |S|}^{-1}} (1 - \pmax(|S|))$, we 
get 
$\P\Bigl(\mbox{$A$ avoids $B$} \,\Bigl|\,
 S \gtom  S' \Bigr) \leq 
\alpha^{-d} {1 + \pmax(|S|) \over 1 - \pmax(|S|)}.$
\end{proof}

\section{Appendix B}
\lab{mlsec}
The purpose of this  section is 
to prove Corollary \ref{cl},
which is used to bound the root profile. 
If $K$ is the transition kernel for 
a Markov chain on the state space $V$, we will consider 
$K$ as an operator acting on the space of functions $f: \Omega \to \R$ by
\be
\lab{op}
Kf(x) = \sum_{y \in \Omega} K(x,y) f(y).
\ee
We will need the following lemma, which was proved
by Yuval Peres.
\begin{lemma}
\lab{yp}
Let $K$ be a doubly stochastic transition kernel and define 
$\kt = K K^t$.
For any function $g: \Omega \to [0,1]$ and $n \geq 1$ we 
have
\[
\ss{K^t g} \leq \la g, g \ra^{1 - {1 \over n}}
\la \kt^{n} g, g \ra^{1 \over n}.
\]
\end{lemma}
\begin{proof}
Since $\kt$ is symmetric it is diagonalizable. Thus
we can write $g = \sum_i \alpha_i g^i$, where the 
$g^i$ are orthonormal eigenfunctions of $\kt$
with corresponding eigenvalues $\lambda_i$. We have
\begin{eqnarray}
{\ss{K^t g} \over \la g, g \ra} &=&
{\la \kt g, g \ra  \over \la g, g \ra} \\
&=& 
{ \sum_i \alpha_i^2 \lambda_i  \over 
\sum_i \alpha_i^2} \\
&\leq& 
\label{48}
\Bigl( { \sum_i \alpha_i^2 \lambda_i^{n}  \over 
\sum_i \alpha_i^2} \Bigr)^{1/n} = 
\Bigl({ \la \kt^{n} g,  g \ra \over \la g, g\ra }
\Bigr)^{1/n},
\end{eqnarray} 
by Jensen's inequality. 
Multiplying both sides by $\la g, g \ra$ yields the lemma. 
\end{proof}
We will also need the following lemma, which was proved by Keith
Ball. 
\begin{lemma}
\label{normlemma}
Let $X$ be a random variable taking values in $[0,1]$ and suppose that
$\e(X) = \mu \le \half$.  
Then for any $p > 1$ we have
\be
\lab{norm}
{\e(X^p) \over \mu^p} - 1 \le (\mu^{1-p} - 1) 
\e\Bigl|{X - \mu \over
  \mu}\Bigr|.
\ee
\end{lemma}
\begin{proof}
Let $l = \half \e(|X - \mu|)$.
For a
given value of $l$, the l.h.s.~of (\ref{norm}) 
is maximized when $X$ is concentrated on
the three values $0, \mu$ and $1$ (because it
is a convex
function of $X$). Let $p_0, p_{\mu}$ and $p_1$ be the respective
probabilities. Then $l = p_1 (1 - \mu) = p_0 \mu$, and 
hence $p_{\mu} = 1 - p_0 - p_1 = 1 -{l \over \mu(1-\mu)}$. It follows
that
\begin{eqnarray*}
{\e(X^p) \over \mu^p} - 1 &=&
{p_1 + p_{\mu} \mu^p \over \mu^p} - 1 \\
&=& l  \Bigl[ {1 \over \mu^p (1 - \mu)} - 
{1 \over \mu (1 - \mu)} \Bigr] \\
&\le& {2l \over  \mu} (\mu^{1-p} - 1),
\end{eqnarray*}
since $1 - \mu \geq \half$, and the proof is complete.
\end{proof} 

Fix $\ds \leq d$. Recall that the $\ds$--truncated Thorp shuffle
is the Markov chain with transition kernel
$\kst = K_1 \dots K_\ds$. 
Let $V$ denote the state space of this chain. 
Corollary \ref{cl} is a consequence of the 
following technical lemma. 
\begin{lemma}
\label{mainlemma}
Fix $f: \Omega \to [0,1]$. 
Then there is a universal constant $C \in (0,1)$ such that 
\[
\ss{\kstarp{t} f} 
\leq \lonep{f}{1 + 1/C d^2 \dsp{12}}.
\]
\end{lemma}
\begin{proof}
Suppose that $\ds=1$. 
Then the truncated Thorp shuffle 
makes the distribution uniform over $V$ in one step. Thus, 
\begin{eqnarray}
\ss{\kstarp{t} f} 
&=& 
\sum_{x \in V} \lonep{f}{2} {1 \over |V|} \\
&=& \lonep{f}{2} \\
\label{basecase} 
&\leq&  \lonep{f}{p},
\end{eqnarray}
for any $p \in [1,2]$, since $\lone{f} \leq 1$.
Suppose now that $\ds \geq 2$.   
Let $c$ be the constant appearing in Corollary \ref{trunccor}.
We will consider the cases $\lone{f} \leq 6^{-c \dsp{6}}$ and
$\lone{f} > 6^{-c \dsp{6}}$ 
separately.\\
\\
{\bf Case 1: $\lone{f} \leq 6^{-c\dsp{6}}$.}
We show by induction on $\ds$ that 
$\ss{\kstarp{t} f} \leq \lonep{f}{1 + 1/c d\dsp{5}}$.
The base case $\ds=1$ is handled by equation (\ref{basecase}) above.
Now assume that the result holds for $\ds-1$. 
Define $\lds$ as the set of vertices in the cube whose $\dsp{{ \, th}}$ 
coordinate is $0$.
Let $\b$ denote the collection  of 
subsets $b$ of $\{1, \dots, 2^d\}$ 
such that 
$X(b) = \lds$ for some $X \in V$ (i.e., there is a configuration 
$X \in V$ such that the set of cards occupying $\lds$ is $b$).
For $b \in \b$, define
$\Omega_b = \{X \in \Omega: X(b) = \lds \}$.
Let $r = \lone{f}$ and
for $\Lambda \subset \b$, define
\[
\Omega_\Lambda = \cup_{b \in \Lambda} \Omega_b.
\]
Let
\[
H = \Bigl\{ b \in \b: { \lone{f  \one_{\Omega_b}}  
\over \lone{f}} \geq {r^{-1 / \ds} \over |\b|} 
\Bigr\}.
\]
Since $
\sum_{b \in \b} {
 \lone{f \one_{\Omega_b}}
\over \lone{f}} = 
 {\lone{f}
\over \lone{f}} = 1$,
Markov's inequality implies that
\be 
\lab{mar1}
{|H| \over |\b|} \leq r^{1/\ds}.
\ee
Let $A = \Omega_H$
and let $f_1 = f \chi_A$ and $f_2 = f \chi_{A^c}$. Then
\be
\lab{sums}
\ss{ \kst f} = \ss{ \kst f_1 + \kst f_2} \leq 2 \ss{\kst f_1}
+
2 \ss{\kst f_2}.
\ee 
We will bound each term on the right hand side separately. First, consider
$\ss{\kst f_1}$. Let $\kt$ be the transition kernel for the
$\ds-$truncated zigzag shuffle,
i.e,  $\kt = K_1 \cdots K_{\ds} \cdots K_1$. 
Let $n = cd (\ds - 1)^4$.
Using Corollary \ref{trunccor} (with $k=1$) and  
combining this with equation (\ref{mar1})
gives
$\kt^{n}(x, \Omega_H) \leq \exp(\exp(-1)) r^{ {1/\ds}}$ 
for all $x$. 
Hence
\begin{eqnarray}
\la \kt^{n} f_1, \chi_{A} \ra &=& {\textstyle |V|^{-1}}
\sum_x f_1(x) \kt^{n}(x, \Omega_H) \\
\lab{ff}
&\le& \lone{f_1} \exp(\exp(-1)) r^{{1/\ds}}.
\end{eqnarray}
Finally, Lemma \ref{yp} gives
\begin{eqnarray}
\ss{ \kst f_1} &\le& 
\la f_1, f_1 \ra^{1 - {1/n}}
\la \kt^{n} f_1, f_1 \ra^{1/n} \\
&\le& 
\la f_1, f_1 \ra^{1 - {1/n}}
\la \kt^{n} f_1, \chi_A \ra^{1/n},
\end{eqnarray}
where the second inequality holds because $f_1 \leq \chi_A$.
Putting this all together, we get
\begin{eqnarray}
\ss{ \kst f_1}  &\leq& \la f_1, f_1 \ra^{1 - {1 /n}} \,
\Bigl[ { 
\lone{f_1} (\exp(\exp(-1)))r^{{1/\ds}} }
\Bigr]^{1/n} \\
&\leq& 
2 \Bigl( { \la f_1, f_1 \ra \over \lone{f_1}} \Bigr)^{1 - {1/ n}}
\times \lone{f_1} \times
r^{{1/\ds n}},
\end{eqnarray}
since $\exp({1 \over n} \exp(-1)) \leq 2$ for all $n$.
Since $n = c d (\ds-1)^4$, and 
${\la f_1, f_1 \ra \over \lone
{f_1}} \leq 1$, we have
\begin{eqnarray}
\label{f1bound} 
\ss{ \kst f_1} 
\le
2 r^{1/cd \ds (\ds-1)^4} \lone{f}.
\end{eqnarray}

Next we bound $\ss{\kst f_2}$. 
Since $K_{\ds}$ is symmetric it contracts $l^2$. Hence
\begin{eqnarray}
\ss{\kst f_2} &\leq& 
\ss{K_{(\ds - 1)} \cdots K_1 f_2} \\
\label{trian}
&=& \sum_{b \in \b} \ss{ K_{(\ds-1)} \cdots K_1 f_2 \one_{\Omega_b}}.
\end{eqnarray}
Note that $K_1 \cdots K_{(\ds - 1)}$ is just the transition kernel for a 
$(\ds-1)$-truncated 
Thorp shuffle and that the $\Omega_b$
are communicating classes for this process. 
Thus, we can use the induction hypothesis to bound each
$\ss{ K_{(\ds-1)} \cdots K_1 f_2 \chi_{\Omega_b}}$, provided that
the corresponding normalized $l_1$ norm
${\lone{f_2 \chi_{\Omega_b}}  \over \lone{\chi_{\Omega_b}}}$
is sufficiently small. 
Define
$r_b := {\lone{f_2 \chi_{\Omega_b}}  \over \lone{\chi_{\Omega_b}}}$.
We claim that for every $b \in \b$ we have
$r_b \leq 
r^{\ds -1 \over \ds}$. 
To see this, note that if  $b \in H$, then
$\lone{f_2 \chi_{\Omega_b}} = 0$ and the claim holds trivially, so assume 
$b \notin H$. Then
\begin{eqnarray}
{\lone{f_2 \chi_{\Omega_b}}  \over \lone{\chi_{\Omega_b}}}
&=&
{\lone{f_2 \chi_{\Omega_b}}} |\b| \\
&\le& \lone{f \chi_{\Omega_b}} |\b| \\
\lab{cool}
&\le& r^{-1 \over \ds} \lone{f}  = r^{\ds-1 \over \ds},
\end{eqnarray}
where the first equality holds because 
$\lone{ \one_{\Omega_b}} = |\b|^{-1}$, 
the second inequality holds because $b \notin H$
(and by the definition of $H$)
and last equality holds because $\lone{f} = r$.
It follows that
\be
\lab{in}
r_b \leq r^{\ds-1 \over \ds} \leq 6^{- c (\ds-1)\ds^5}
\leq 6^{-c (\ds-1)^6}. 
\ee
Thus we can apply the
induction hypothesis, which gives
\begin{eqnarray}
\ss{ K_{(d'-1)} \cdots K_1 f_2 \one_{\Omega_b}}
&\leq& 
r_b^{ 1/c d(\ds-1)^5} 
\lone{f_2 \chi_{\Omega_b}} \\
&\le& r^{ 1/cd \ds (\ds-1)^4} 
\lone{f_2 \chi_{\Omega_b}},
\end{eqnarray}
where the second inequality 
follows from the first inequality in 
(\ref{in}).  
Combining this with equation (\ref{trian}) 
and using the fact that $f_2 \leq f$
gives
\be
\lab{f2bound}
\ss{\kst f_2} 
\leq  r^{1/cd \ds (\ds-1)^4}
\lone{f}.
\ee
We are now ready to bound $\ss{ K^t f}$. 
Combining equations (\ref{f2bound}), 
(\ref{f1bound}) and (\ref{sums}),  we get
\begin{eqnarray} 
\ss{ \kst f} &\leq& 
\label{fone}
\Bigl( 6 r^{1/cd \ds (\ds-1)^4} \Bigr) \lone{f}.
\end{eqnarray}
Since $(k-1)^{-4} - k^{-4} \geq k^{-5}$
for integers $k \geq 2$,
the quantity (\ref{fone}) is at most
\[
 6 r^{1/cd \dsp{5} + 1/c d \dsp{6}} \lone{f} \leq
r^{1/cd \dsp{5}}
\lone{f},
\]
since  $r \leq 6^{-cd \dsp{6}}$.
This concludes the proof in the case 
$r \leq 6^{-cd \dsp{6}}$.
\\
{\bf Case 2: 
$r > 6^{-cd \dsp{6}}$.}
Let $C$ be an integer that is larger than 
$2^{15} c^2 15 \log 2 \log 6$.
We will show by induction on $\ds$ that
\[
\ss{\kst f} \leq r^{1 + 1/C d^2 \dsp{12}}.
\]
The base case $\ds =  1$ was handled earlier by equation (\ref{basecase}).

Now, fix $\ds \geq 2$ and $f: \Omega \to [0,1]$ and suppose that
$r = {\lone{f}} 
> 6^{-cd \dsp{6}}$.
We
can assume w.l.o.g.~that $r \leq \half$.
Otherwise, let $h = 1 - f$, and suppose that the result holds for $h$, i.e.,
for $q = 1/C d^2 \dsp{12}$ we have
\[
\ss{\kst h} \leq \lonep{h}{1 + q},
\]
or equivalently,
\be
\lab{estat}
{ \lone{h} - \ss{\kst h} } \geq 
\Bigl[ 1 -  {\lonep{h}{q} }
 \Bigr] {\lone{h} }.
\ee
Note that 
\begin{eqnarray}
\ss{ \kst h} &=& \la 
\kst (1 - f), 
\kst (1 - f) \ra \\
&=& \la \kst \chi, \kst \chi \ra 
- 2 \la \kst \chi, \kst f \ra + \la \kst f, \kst f \ra \\
&=& 1 - 2 \lone{f} + \ss{\kst f}\\
 &=& \lone{h} - \lone{f} + \ss{\kst f},
\end{eqnarray}
where the third equality holds because $\kst$ is doubly stochastic and hence
$\kst \chi = \chi$. 
Thus
\be
\lab{same}
\lone{h} - \ss{K^th} 
=
\lone{f} - \ss{K^tf}.
\ee
Define $u:[0,1] \to \R$ by
\be
\lab{udef}
u(x) = (1 - x^{q}) x = { x(1 - x) \over 1 + x^q + 
\cdots + x^{1-q}},
\ee
so the RHS of (\ref{estat}) is $u( {\lone{h}})$.
Since the 
numerator on the RHS of (\ref{udef}) 
is symmetric about $\half$ and the denominator is 
increasing, we have $u(x) \geq u(1-x)$ if $x \leq \half$. 
This, combined with equation (\ref{same}), shows that
equation (\ref{estat}) is still true if 
we replace the $h$ by $f$. Thus we can assume henceforth that
$r \leq \half$. 
 
Let $\b$ and $\Omega_b$ 
be as defined above.
Then
\begin{eqnarray}
\ss{ \kst f } 
&\le& \ss{ K_{(\ds -1)}  \cdots K_1 f } \\
\lab{tophalf}
&=& \sum_{b \in \b} \ss{ K_{(\ds - 1)} \cdots K_1 f \chi_{\Omega_b} }.
\end{eqnarray}
For $b \in \b$, define $r_b = {\lone{ f\chi_{\Omega_b}} \over 
\lone{\chi_{\Omega_b}}} =
{\lone{ f\chi_{\Omega_b}}} \,|\b| $. 
We may assume that
\[
\ss{ K_{(\ds - 1)} \cdots K_1 f \chi_{\Omega_b} }
\leq r_b^{1/C d^2 \dsmp{12}} \lone{f \one_{\Omega_b}}.
\]
(In the case where 
$r_b \leq 6^{-cd \dsmp{6}}$
this was proved earlier,
since $c d \dsmp{5} \leq C d^2 \dsmp{12}$; in the case where
$r_b > 6^{-cd \dsmp{6}}$
this is the induction hypothesis.)
Combining this with (\ref{tophalf}) gives
\[
\ss{ \kst f } 
\le  \sum_{b \in \b} 
r_b^{1/C d^2 \dsmp{12}}
\lone{f \one_{\Omega_b}}
= |\b|^{-1}
\sum_{b \in \b} 
r_b^{1 + 1/C d^2 \dsmp{12}}.
\]
Thus, unless
\be
\lab{ass}
|\b|^{-1}
\sum_{b \in \b}
r_b^{1 + 1/C d^2 \dsmp{12}} 
\geq
r^{1 + 1/C d^2 \dsp{12}},
\ee
the result is immediate. So assume that 
(\ref{ass}) holds. 
For $b \in \b$, define $w_b = {\lone{f \chi_{\Omega_b}} \over \lone{f}}$.
Note that $\sum_{b \in \b} w_b = 1$. 
Let $U$ be chosen  uniformly at random  from $\b$.
Let $p = 1 + 1/Cd^2\dsmp{12}$.
Dividing both sides 
of (\ref{ass}) by $r^{p}$ gives
\begin{eqnarray}
\lab{ttt}
{\e(r_U^{p}) 
\over r^{p}
} \geq r^{1/Cd^2 \dsp{12} -1/Cd^2 \dsmp{12}}.
\end{eqnarray}
Using the inequality
$k^{-12} - (k-1)^{-12} \leq - k^{-13}$, valid 
for integers $k \geq 2$, and subtracting $1$ from both sides of
(\ref{ttt}) 
gives
\be
\lab{forest}
{\e(r_U^{p}) \over 
r^{p}} - 1 \geq r^{-1/Cd^2 \dsp{13}} - 1
\ee
Let $\u$ be the uniform probability measure on $\b$
and let $\newmu$ be the measure on $\b$ defined by the $w_b$.
Define
\[
\Vert \u - \newmu \Vert_{TV} = \half \sum_{b \in \b} 
\Bigl| w_b - {\textstyle |\b|^{-1}} \Bigr| = \half \e\Bigl| {r_U - r \over r}
  \Bigr|.
\]
Note that
$\e(r_U) = |\b|^{-1} \sum_{b \in \b} r_b = r$. 
Plugging $X = r_U$ and $\mu = r$ into 
Lemma \ref{normlemma} and combining with equation (\ref{forest}) gives
\begin{eqnarray}
2 \Vert \u - \newmu \Vert_{TV} &\geq& {
r^{-1/C d^2 \dsp{13}} - 1 \over
r^{-1/C d^2 \dsmp{12}} - 1} \\
&=& 
\label{qexp}
{\exp \Bigl( { - \log r \over C d^2 \dsp{13}} \Bigr) - 1 \over 
\exp \Bigl( { - \log r \over C d^2 \dsmp{12}} \Bigr) - 1}.
\end{eqnarray}
Since 
$r > 6^{-cd \dsp{6}}$, the quantities in the exponents in 
(\ref{qexp}) are in $(0,\half]$. 
(Recall that $C$ is much larger than $c$.)
Hence, the fact that
${e^t - 1 \over t} \in [1,2]$ whenever $t \in (0,\half]$
implies that the quantity in (\ref{qexp}) is at least
\[
{\dsmp{12} \over 2 \dsp{13}} \geq {\dsp{-1} \over 2^{13}},
\]
where the inequality holds because $\ds \geq 2$ and hence
${\ds - 1 \over \ds} \geq \half$.
It follows that $\Vert \u - \newmu \Vert_{TV} 
\geq 2^{-14} \dsp{-1}$.
Note that 
\begin{eqnarray}
2 \Vert \u - \newmu \Vert_{TV} &=& 
\sum_{b \in \b} 
\max(\newmu(b), \u(b)) -  
\min(\newmu(b), \u(b)); \\
2 &=&   
\sum_{b \in \b} 
\max(\newmu(b), \u(b)) +  
\min(\newmu(b), \u(b)).
\end{eqnarray}
Subtracting the first equation from the second and dividing by $2$ gives
\be
\lab{tvsum}
1 -   \Vert \u - \newmu \Vert_{TV}
= \sum_{b \in \b} \min(\newmu(b), \u(b)).
\ee 
Recall that $\kt$ is the transition kernel for the $\ds$--truncated zigzag
shuffle. 
Note that 
\begin{eqnarray}
\la f, 
\kt^{n} f \ra &=& 
\sum_{b \in \b} \la f \one_{\Omega_b}, (\kt^{n} f ) \one_{\Omega_b}
\ra \\
&\leq& \sum_{b \in \b} \min\Bigl( \lone{f \chi_{\Omega_b}}, 
\lone{(\kt^{n} 
f) \chi_{\Omega_b}} 
\Bigr) \\ 
&=& \lone{f} 
\sum_{b \in \b} \min\Bigl( w_b, {\lone{(\kt^{n} f) \chi_{\Omega_b}} 
\over \lone{f}} 
\Bigr),
\end{eqnarray}
where the inequality holds because 
$f \chi_{\Omega_b} \leq 1$ and  $(\kt^{n} f) \chi_{\Omega_b} \leq 1$. 
Let $n = 15 c d \dsp{5} \log 2$.
Corollary  \ref{trunccor} implies that 
${ \lone{ (\kt^{n} f) \chi_{\Omega_b}} 
\over \lone{f}}  \leq \alpha(\ds) |\b|^{-1}$,
where $\alpha(k) := \exp(2^{-15k})$. 
Hence, 
\begin{eqnarray}
\la f, 
\kt^{n} f \ra &\leq& \lone{f} \,\et 
\sum_{b \in \b} \min( w_b, \sfrac{1}{|\b|}) \\
&=& \lone{f} \,\et (1 - \Vert \newmu - \u \Vert_{TV}) \\
&\leq& \lone{f} \,\et \left[1 - 2^{-14} \dsp{-1} \right].
\end{eqnarray}
Hence Lemma \ref{yp} 
gives
\begin{eqnarray}
\ss{ \kst f} &\leq& \la f, f \ra^{1 - {1/ n}} \,
\Bigl[ {  
\lone{f}\, \et (1 - {2^{-14} \dsp{-1}} )}
\Bigr]^{1 /n} \\
&=& 
\Bigl( { \la f, f \ra \over \lone{f}} \Bigr)^{1 - {1/ n}}
\times \lone{f} \times
\et^{1/n}
\times
(1 - {2^{-14} \dsp{-1}})^{1 /n}
\\
\label{first}
&\le&
\lone{f} \exp \Bigl( {1 \over n} \Bigl[ 2^{-15\ds} - \dsp{-1}{2^{-14}} 
\Bigr] \Bigr) \\
&\le&
\lone{f}
\exp \Bigl( {-1/ 2^{15} c d \dsp{6}15 \log 2} \Bigr),
\end{eqnarray}
since ${\la f, f\ra \over  \lone{f}} \leq 1$
and $2^{-15k} \leq 2^{-15}k^{-1}$ for all positive integers $k$.
Finally, since $r > 6^{-cd \dsp{6}} = \exp(-c d \dsp{6} \log 6)$,
we have $r^{1/C d^2 \dsp{12}} \geq 
\exp( {-1/ 2^{15} c d \dsp{6}15 \log 2}).$ (Recall 
that $C > 2^{15} c^2 15 \log 2 \log 6$.)
It follows that $\ss{\kst f } \leq r^{1/Cd^2 \dsp{12}} \lone{f}$. 
This completes the proof.
\end{proof}
To bound the root profile, we actually used the following corollary.
\begin{corollary}
\lab{cl}
Fix $S \subset \Omega$ and let
\[
x = {|S| \over (2^d)!}.
\]
Let $\{p(x,y)\}$ be the transition probabilities for a 
round of the Thorp shuffle.
Then there is a universal constant $C> 0$ such that 
\[
\ss{p(S, \cdot)} \leq x^{{1 + C/d^{14} }}.
\]
\end{corollary}
\begin{proof}
Let $f = \chi_S$ 
and 
$\ds = d$ and 
apply Lemma \ref{mainlemma}.
(Note that if $K$ is the transition kernel 
for a round of the Thorp shuffle then $p(S, \cdot) = K^t f$.) 

\end{proof}

\noindent {\bf Acknowledgments.} I am grateful to
K.~Ball, T.~Coulhon, E.~Mossel, C.~Nair, Y.~Peres, A.~Sinclair,
D.~Wilson, P.~Winkler and J.~Zuniga
for invaluable discussions. 

 I want to thank Jessica Zuniga for pointing out an error in the 
conference version of this paper.
Yuval Peres proved Lemma \ref{yp}
and Keith Ball proved Lemma \ref{normlemma}.

  I also want to thank Christian Borgs and Jennifer Chayes 
for giving me the opportunity to spend the year at the Theory Group
of Microsoft Research,
where I did much of this research.


\end{document}